\documentclass{ifacconf}

\usepackage{graphicx}      % include this line if your document contains figures
\usepackage{natbib}        % required for bibliography
\usepackage{color}
\usepackage{amsmath,amsfonts,amssymb, bm, bbm}
\usepackage{latexsym}
\usepackage{tikz-network}

\newcommand{\mc}{\mathcal}
\newcommand{\1}{\mathbbm{1}}
\newcommand{\R}{\mathbb{R}}

\newcommand{\E}{\mathbb{E}}

\def\beq{\begin{equation}}
\def\eeq{\end{equation}}
\newtheorem{theorem}{Theorem}

\newtheorem{example}{Example}%}
\theoremstyle{remark}
\newtheorem{remark}[theorem]{Remark}

\DeclareMathOperator{\Var}{Var}

\graphicspath{ {./data_plot_graph/} }

\begin{document}
\begin{frontmatter}

\title{Systemic risk and network intervention} 

	\thanks[footnoteinfo]{Giacomo Como is also with the Department of Automatic Control, Lund University, Sweden. This work was partially supported by MIUR grant
	Dipartimenti di Eccellenza 2018--2022 [CUP: E11G18000350001], the Swedish Research Council, and by the Compagnia di San Paolo.}
	
	\author[First]{Luca Damonte} 
	\author[First]{Giacomo Como} 
	\author[First]{Fabio Fagnani}

	\address[First]{{Department of Mathematical Sciences, Politecnico di Torino}\\ Corso Duca degli Abruzzi 24, 10129 Torino, Italy\\ 
 (e-mail: \{luca.damonte,giacomo.como,fabio.fagnani\}@polito.it).}

\begin{abstract}                % Abstract of not more than 250 words.
We consider a novel adversarial shock/protection problem for a class of network equilibria models emerging from a variety of different fields as
continuous network games, production networks, opinion dynamic models. The problem is casted into a min-max problem and analytically solved for 
two particular cases of aggregate performances: the mean square of the equilibrium or of its arithmetic mean. The main result is on the shape of
the solutions, typically exhibiting a waterfilling type structure with the optimal protection concentrated in a proper subset of the nodes, depending
significantly on the aggregate performance considered. The relation of the optimal protection with the Bonacich centrality is also considered.

\end{abstract}


\end{frontmatter}
%===============================================================================

\section{Introduction}
Many different network models lead to equilibria configurations describable in the form

\beq\label{equilibrium}
 x=(I-\Lambda)^{-1}Dc
\eeq
where $\Lambda$ is an $n\times n$ sub-stochastic matrix (with spectral radius less then $1$), $D$ a an $n\times n$ diagonal non negative matrix,
$c$ a vector of inputs that can be exposed to external shocks, and the vector $x$ is the system equilibrium configuration.

Such relations appear in several socio-economic models, as the Nash equilibria of an underlying game, or
as the asymptotic configuration of a network dynamics model. Notable examples is that of continuous network games  with linear best reply
(\cite{jackson2015games,bramoulle2014strategic,bramoulle2015games}),
production networks (\cite{carvalho2014micro,acemoglu2010cascades,acemoglu2015networks}), coordination games \cite{},
and opinion dynamics models as the Friedkin-Johnsen model (\cite{friedkin1990social}).
In the economic interpretations, the components of the vector $c$ can have the meaning of marginal benefits of the economic agents
or rather of productivity indices of firms. In the sociological models, instead, they typically represent the initial opinion of the agents. 

A large attention has been recently devoted, particularly in the economic literature \cite{}, to the effects that a perturbation in the vector $c$
can have on the network equilibrium, in particular how shocks at the level of single agents can possibly be amplified by the network interaction,
and propagate to the other agents. These studies have shown the role of the network topology (given by the pattern of non-zero elements of the matrix $\Lambda$)
in determining the extent of these contagion phenomena. The performance typically considered in the literature is the mean squared error for $x$ and 
for the aggregated function $n^{-1}\sum_ix_i$, the arithmetic mean of the equilibrium configuration.

In this paper, we take a further step in this direction and we analyze a more complex model, where shocks are complementary paired with protections, and we
cast it into an adversarial min-max problem.
Specifically, we assume that the vector $c$ has the following structure: $$c_i=\bar c_i+q_i^{-1}\eta_i$$ where $\bar c_i$ is a reference value,
$\eta_i$ is a random variable modeling the shock and $q_i\geq 1$ is the protection actuated on node $i$.

In this paper, we consider min-max optimization problems formulated in terms of the mean squared error of either $x$ or the arithmetic mean
$n^{-1}\sum_ix_i$ and we solve it under an (different) assumption on the correlations of the shocks. Solutions to these optimization problems in $q$
will typically exhibit a 'waterfilling' structure with the optimal solution $q$ concentrated on a limited number of nodes. 
The main message coming from our analysis is that the nodes on which protection has to be taken to minimize the effect of external
shocks is not only function of the network topology but also depends on the correlation patterns of the shocks, as well on the type of functional we are considering.

We conclude presenting a brief outline of the rest of this paper. In Section \ref{model} we formally introduce the min-max problems and in section
\ref{application} we review some motivating applications deriving from network games theory, economic network theory and opinion dynamics. 
In Section \ref{independent_shocks} we study the aforementioned problems in the particular case of independent shocks.
In the last part we conclude with possible developments.

\section{The shock-protection optimization problem}\label{model}
We consider the model described by (\ref{equilibrium}) and we rewrite the assumption made on $c$ in the following compact form:
\beq c=\bar c+Q^{-1}\eta\eeq
where, precisely,
\begin{itemize}
\item $\bar c$ is a given reference vector;
\item $\eta$ is a random vector modeling the shock, having mean $0$ and covariance matrix $\omega\in\Omega$;
\item $Q=\emph{diag}(q)$ where $q\in\mc Q$ is the \emph{protection} vector.
%, characterized by $q_i\geq 1$ for every $i$. 
\end{itemize}
We formulate the following min-max problems:
\beq\label{varx}
 \min\limits_{q\in\mc Q}\max\limits_{\omega\in\Omega}\sum\limits_i \Var[x_i] 
\eeq
\beq\label{varm}
 \min\limits_{q\in\mc Q}\max\limits_{\omega\in\Omega}\Var[n^{-1}\1'x] 
\eeq

$\Omega$ and $\mc Q$ model the way the two adversarial can intervene in the system. For the sake of this paper, we have assumed that shocks are
uncorrelated of total variance $1$: \beq\label{uncor}\Omega=\{\omega=\emph{diag}(\sigma_1^2,...,\sigma_n^2)\;\sum_i\sigma_i^2=1\}\eeq
The set $\mc Q$ is instead assumed to have the form
\beq\label{prot}
\mc Q=\{q\in\R_+^n\;|\; 1\leq q_i \,\forall i, \; \|q\|_2\leq C\}
\eeq
The normalization of the power of the shocks to $1$ is not restrictive as a different threshold can always be absorbed modifying the maximum protection norm $C$.

\section{Motivating examples}\label{application}

We here present some multi-agent models described by the equilibrium relation (\ref{equilibrium}).

\subsection{Continuous coordination games and opinion dynamics models}
We are given a set of players $\mc V$ whose utilities are given by
$$u_i(x)=-\frac{1}{2}\left[\sum\limits_j W_{ij} (x_i-x_j)^2+\rho_i (x_i-c_i)^2\right]$$
where the coefficients $W_{ij}$ and $\rho_i$ are all non negative. This can be interpreted as a continuous variable coordination
game with the presence of anchors $c_i$'s.

Assume that $w_i=\sum_k W_{ik}>0$ for every $i$ and put
$$\Lambda_{ij}=\frac{W_{ij}}{w_i+\rho_i},\; D_{ii}=\frac{\rho_i}{w_i+\rho_i}$$
If the set of nodes $\{i\in\mc V\,|\, c_i>0\}$ is globally reachable, it can be shown that $\Lambda$ has spectral radius less
than $1$ and that the game has just one Nash equilibrium given by formula \eqref{equilibrium}
with $\Lambda$ and $D$ given above and $c$ equal to the vector of anchors.

Related to this, one can consider the dynamical system
$$x(k+1) = \Lambda x(k) + Dc$$
Notice how $x(k)$ converges to the equilibrium expressed in \eqref{equilibrium} under the same assumptions above.
With the extra assumption that $c=x(0)$ the initial opinion vector, this is the celebrated Friedkin-Johnsen model in
opinion dynamics. In this setting, $D_{ii}$ can be interpreted as a measure of how agent $i$ is attached to
its own original opinion.

\subsection{Quadratic games}
A variant of the games considered above are the quadratic games where utilities have the form
\beq\label{quadratic}
 u_i(x)=c_i x_i-\frac{1}{2}x_i^2+ \beta\sum\limits_j W_{ij} x_ix_j
\eeq
where the $c_i$'s and $\beta$ are positive constants and the elements $W_{ij}$ are non negative. They are used \cite{}
to model socio-economic systems with complementarity effects (e.g. education, crime).
The first two terms of (\ref{quadratic}) give the benefits and
the costs to player $i$ of providing the action level $x_i$. The last term instead
reflects the vantage of cooperation of $i$ with his friends (those $j$ for which $W_{ij}>0$). In the case when $\beta w_i<1$
for all $i$ (where, as above, $w_i=\sum_j W_{ij}$), 
we have that the unique Nash equilibrium is given by
$x=(I-\beta W)^{-1}c$, a special case of (\ref{equilibrium}).

\subsection{Production networks}\label{production}
In the Cobb-Douglas model of an economy with firms producing distinct goods interconnected in the production process,
the production outputs of the firms satisfy at equilibrium the following relation
\beq\label{CD}
 x_i=\beta\sum P_{ij} x_j+(1-\beta) c_i
\eeq
Here $x_i$ has to be interpreted as the log of the output of firm $i$, $P_{ij}$ as the share of the good produced by
$j$ in the production technology of firm $i$, $c_i$ is the log productivity shock to firm $i$, and finally $\beta<1$
is a constant indicating the level of interconnection in the economy. Relation (\ref{CD}) can be rewritten as $x=(I-\beta P)^{-1}(1-\beta) c$.
In this model, a natural candidate for the economy's performance that is the real value added, can be expressed (in logarithm)
simply as $y=n^{-1}\sum_i x_i$.

\section{Main results}\label{independent_shocks}
In this section, we give a complete solution to the optimization problems (\ref{varx}) and (\ref{varm}) under the standing assumption
previously made of uncorrelated shocks, assuming that $\Omega$ and $\mc Q$ have the form (\ref{uncor}) and (\ref{prot}), respectively.

In this case, the two cost functions can be rewritten as follows. We put $L=(I-\Lambda)^{-1}D$ and we compute
\beq\label{varx2}
 \begin{array}{rcl}\sum\limits_i \Var[x_i] &=&\E[\eta'Q^{-1}L'LQ^{-1}\eta]\\
  &=&{\rm tr}(Q^{-1}\E[\eta\eta']Q^{-1}L'L)\\ &=&\sum\limits_i \left( \sigma_i \frac{\ell_i}{q_i} \right)^2
 \end{array}
\eeq
where $\ell_i=\sqrt{(L'L)_{ii}}$ is the squared $2$-norm of the $i$-th column of the matrix $L$.

Similarly, if  we put $v=n^{-1}L'\1$, the variance of the arithmetic mean can be computed as
\beq\label{varm2}
 \Var[n^{-1}\1'x]=v'Q^{-1}\E[\eta\eta']Q^{-1}v=\sum\limits_i \left( \sigma_i \frac{v_i}{q_i} \right)^2
\eeq

The two min-max problems can thus be framed under the general problem
\beq\label{indi_varx}
 \min\limits_{\|q\|\leq C}\,\max\limits_{\|\sigma\|\leq 1}\, \sum\limits_i \left( \sigma_i \frac{y_i}{q_i} \right)^2,
\eeq
where $y=\ell$ in the case of problem (\ref{varx}) and $y=v$ in the case of problem (\ref{varm}).
We will refer to $y_i$ as to the \emph{centrality} of agent $i$ and, without lack of generality, we assume that the elements of $y$ are ordered in a decreasing order,
i.e $y_1 \geq y_2 \geq ... \geq y_n$.

We now solve problem (\ref{indi_varx}). First we introduce the function $f:(0,+\infty)\to\R$ given by
\beq\label{function}f(\lambda)= \sum_{i=1}^n \mbox{ max}\left\{ 1, \left( y_i / \sqrt{\lambda} \right)^2\right\}.\eeq
We notice that $f$ is continuous, strictly decreasing in $ (0, y_1^2]$, and 
$$\lim\limits_{\lambda\to 0+} f(\lambda)=+\infty,\quad f(y_1^2)=n.$$
This implies that for every $C\geq\sqrt n$, it is well defined $\lambda(C):=f^{-1}(C^2)$.
We now define $k(C)$ as the maximum index such that $ y_{k(C)}>\sqrt{\lambda(C)} $.

The following result holds true
\begin{theorem}\label{prop_indi} Let $C\geq\sqrt n$.
 It holds
 \beq\label{indi_varx_prop}
  \min\limits_{\|q\|\leq C}\,\max\limits_{\|\sigma\|\leq 1}\, \sum\limits_i \left( \sigma_i \frac{y_i}{q_i} \right)^2=\lambda(C)
 \eeq
  and the optimum value for $q$ is reached by
 \beq\label{waterfilling_indi_varx} 
  q_k = \left\{\begin{array}{ll}y_k /\sqrt{\lambda(C)}\quad &{\rm if}\, k\leq k(C)\\ 1 \quad &{\rm otherwise}
  \end{array}\right.
 \eeq
\end{theorem}

\begin{pf}

Notice first that the internal maximization is solved by any vector $\sigma$ concentrated on nodes $i$ for which $y_i/q_i$ is maximal.
We can thus reformulate the problem as
\beq\label{indi_varx2}
 \min\limits_{\|q\|\leq C}\max\limits_i \left( y_i/q_i \right)^2.
\eeq
and, introducing a new slack variable $\psi \in \R$, finally as
 \begin{align*}
  \min\limits_{\psi, \|q\|\leq C} \,& \psi\\
    \emph{s.t.} \,& \left( y_i/q_i \right)^2 \leq \psi, \, i=1,...,n. 
 \end{align*}
 We introduce Lagrange multipliers $\alpha \in \R^n$, $\delta \in \R^n$, and $\gamma$ and the Lagrangian function
 \begin{align*} L(q,\psi,\alpha,\delta,\gamma) = \psi +& \sum_{i=1}^n \alpha_i \left[ \left( y_i/q_i \right)^2 - \psi \right] +\\ 
                                       &\gamma \left( \sum_{i=1}^n q_i^2 - C^2 \right) - \sum_{i=1}^n \delta_i \left( q_i - 1 \right).
 \end{align*}
 Since the objective function and the constraints are convex and differentiable Karush-Khun-Tucker (KKT) conditions
 are necessary and sufficient for finding the optimum.
 Indicating with $(q^*,\psi^*)$ and $(\alpha^*,\delta^*,\gamma^*)$ the optimum points for, respectively, the primal and the dual problem,  
 KKT conditions are expressed as
  \begin{align}\label{KKT1}
   q_i^* \geq 1,\quad \alpha_i^* \geq 0,\quad \left( y_i/q_i^* \right)^2 &\leq \psi^*, &i=1,...,n \\
   \label{KKT2}
   \left( (q_i^*)^2 - y_i \sqrt{\beta_i^*/\gamma^*} \right) \left( q_i^*-1 \right) &= 0, &i=1,...,n \\
   \label{KKT3}
   \beta_i^* \left[ \left( y_i/q_i^* \right)^2 - \psi^* \right] &= 0, &i=1,...,n \\
    \label{KKT4}
   \| q^* \|_2 &= C,\\
   (q_i^*)^2 - y_i \sqrt{\alpha_i^*/\gamma^*} &\geq 0, &i=1,...,n \\
   1 - \sum_{i=1}^n \alpha_i^* &= 0, &i=1,...,n .
 \end{align}
 Note that, if for some $i$,  $\psi^*<y_i^2$, then, by (\ref{KKT1}),  $q_i^*>1$. From (\ref{KKT2}) we deduce that $\beta_i^*\neq 0$
 and (\ref{KKT3}) finally yields $q_i^*=y_i^2/\psi^*$.
 If instead $\psi^*\geq y_i^2$, then $q_i^*=1$. In fact, if $q_i^*>1$ then, from (\ref{KKT3}) it would follow that $\beta_i^*= 0$, and, substituting in (\ref{KKT2}), we would get a contradiction.
 
 Therefore, for every $i$,  $q_i^*=\max\left\{ 1, y_i / \sqrt{\psi^*} \right\}$. Plugging these values in (\ref{KKT3})
 we obtain $f(\psi^*)=C^2$ or, equivalently,

 $\psi^* = \lambda(C)$.
 This proves the result.
 \qed
\end{pf}

We now comment on the result obtained.

\begin{remark}
The structure of the optimum we have found shows how the optimal protection is in general concentrated on a proper subset of nodes. In this respect,
it is interesting to analyze various regimes depending on the chosen budget cost $C$.
\begin{enumerate}
\item Notice that protection is active in just one node, namely $q_k=1$ for all $k>2$, if and only if $y_2<\sqrt{\lambda(C)}$ or, equivalently, 
$$C<\sqrt{f(y_2^2)}= \sqrt{n+y_1^2/y_2^2-1}$$
In this case, we get from (\ref{function}) that the optimal value is given by
\beq\label{lowbudget}\lambda(C)=y_1^2/(C^2-n+1)\eeq
We will refer to this as to the \emph{low budget regime}.
\item Notice that protection is active on all nodes, namely $q_k>1$ for all $k$, if and only if  $y_n>\sqrt{\lambda(C)}$ or, equivalently, 
\beq\label{threshold}C>\sqrt{f(y_n^2)}= \|y\|/y_n\eeq
We will refer to this as to the \emph{high budget regime}.
In this regime, we get from (\ref{function}) that the optimal value si given by
\beq\label{?highbudget}\lambda(C)=\|y\|^2/C^2\eeq
\item It follows from the shape of the optimum (\ref{waterfilling_indi_varx} ) that when the protection is active $q_i>1$, then the level of protection $q_i$
is proportional to the centrality of the node $y_i$. 
\end{enumerate}

\end{remark}

\section{Examples}
We now go back to our original problems and present a number of examples. Clearly, the optimal protection solutions to the two optimization problems (\ref{varx2})
and (\ref{varm2}) in the measure that the two centrality vectors, either $\ell$ for the minimization of the mean square error of $x$ or $v$ for the minimization
of the mean square error of $n^{-1}\1'x$, may differ. In the following we denote by $q_v^*$ and $q_{\ell}^*$ the optimal protection solutions related to the two problems.

For the sake of simplicity we  consider a special case of our general network equilibrium problem, already met in the motivating examples:
\beq\label{equilibriumspecial}
 x=(1-\beta)(I-\beta P)^{-1}c
\eeq
where $P$ is a stochastic matrix and $\beta\in\ (0,1)$. It corresponds to the model for the production output of a networked economy reviewed in Subsection \ref{production}.
It can also be seen as the Nash equilibrium of a coordination game in the special case when the strength of the anchors relatively to the network interactions strength is the same for all nodes.

In this case, particularly in the economic applications,  the matrix $L=(1-\beta)(I-\beta P)^{-1}$ (that is also stochastic) takes the name of \emph{Leontief} matrix.
The vector $v=n^{-1}L'\1$ is a stochastic vector, known as the {\emph{Bonacich centrality} of the network described by $P$.

For this problem, we study the behavior of $q_v^*$ and $q_{\ell}^*$ for three different networks: a star undirected graph, a directed graph with $11$ nodes,
and a real dataset of a production network with $417$ node.

\begin{example}
  
Consider the undirected star graph $S_{n+1}$ with $n+1$ nodes. Denoted by $A$ its adjacency matrix, we put $P_{ij}=d_i^{-1}A_{ij}$
where $d_i=\sum_jA_{ij}$ is the degree of node $i$.
Notice that for symmetric reasons, all quantities ($v_i$, $\ell_i$, $(q_v^*)_i$, $(q_{\ell}^*)_i$) will be the same for all leaves $i$ in the star.
We use subscripts $0$ and $\varepsilon$ to indicate any component relative, respectively, to the center of the star or to the leaves.

  Values of $v$ and $\ell$ are:
  \begin{align*}
     v_0 &= \frac{1+\beta(n-1)}{n(1+\beta)},\quad v_{\varepsilon} = \frac{\beta+(n-1)}{n(n-1)(1+\beta)} \\
     \ell_0 &= \left( \frac{1+\beta^2(n-1)}{(1+\beta)^2} \right)^{1/2}\\
     \ell_{\varepsilon} &= \left( \frac{\beta^2+(n-1)\left[ 2(2\beta^2-\beta^4)+n(1-\beta^2)^2-1\right]}{(n-1)^2(1+\beta)^2} \right)^{1/2}.
  \end{align*}
For such family of graphs, there will be just one threshold describing the behavior of each optimal solution as shown in  (\ref{threshold}):
$\|v\|/v_{\epsilon}$ and $\|\ell\|/\ell_{\epsilon}$. When $C$ is below this threshold, the corresponding $q^*$ (either $q^*_v$ or $q^*_{\ell}$)
is concentrated in just the center of the star, when is above, it is instead diffused also on the leaves. Computing such thresholds for large $n$, we obtain that
$$\|v\|/v_{\epsilon}\sim\beta n,\quad  \|\ell\|/\ell_{\epsilon}\sim\gamma \sqrt n$$
where $$\gamma=\frac{(1+2\beta^2-2\beta)^{1/2}}{1+\beta}$$
This significative difference between the two thresholds indicates that there are going to be three different regimes for the budget with a central
large regime for which the optimal protection for the mean square error of $x$ is diffused, while the one for the mean square error of $n^{-1}\1'x$
is instead concentrated on the center of the star. This is intuitive as in the aggregate performance $n^{-1}\1'x$, shocks at leaves level have a smaller
impact because of the cancellation induced by averaging.
In table \ref{tb:table1} we report, in the large scale limit, the details of the protection vectors in the three regimes.

  \begin{table}[!hb]
  \begin{center}
  \caption{Optimal protections $q_v^*$ and $q_{\ell}^*$ for star graph when $n\rightarrow \infty$.}  \label{tb:table1}
   \begin{tabular}{|l|l|l|l|}
    \hline
     & $C \leq \gamma \sqrt{n}$ & $C \in \left(\gamma \sqrt{n},\beta n\right)$ & $C \geq \beta n$ \\ %\frac{\beta^2}{(1-\beta^2)^2}
    \hline\hline
    $q_{v_0}^*$ & $\sqrt{C^2-n}$ & $\sqrt{C^2-n}$ & $C$ \\ \hline
    $q_{v_{\varepsilon}}^*$ & 1 & 1 & $C/(\beta n)$ \\ \hline
    $q_{\ell_0}^*$ & $\sqrt{C^2-n}$ & $C \sqrt{1-1/\gamma^2}$  & $C \sqrt{1-1/\gamma^2}$ \\ \hline
    $q_{\ell_{\varepsilon}}^*$ & 1 & $C/\left(\gamma \sqrt{n}\right)$ & $C/\left(\gamma \sqrt{n}\right)$ \\ \hline 
   \end{tabular}
  \end{center}
 \end{table}
 
\end{example}

\begin{example}

 Consider the directed network of $n=11$ agents in figure \ref{fig:network}.
 
 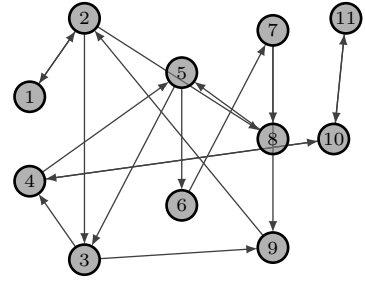
\begin{figure}[!ht]
  \begin{center}
   \begin{tikzpicture}[scale=0.8]
    \Vertices[size=0.4,opacity=0.6,color=gray]{./data_plot_graph/vertices.csv}
    \Edges[lw=0.5,NoLabel]{./data_plot_graph/edges.csv}
   \end{tikzpicture}
   \caption{Directed network of $n=11$ nodes.}
   \label{fig:network}
  \end{center}
 \end{figure}
 As in the previous example we take $P$ as the normalized adjacency matrix $P_{ij}=d_i^{-1}A_{ij}$. 
 In this example we have considered $\beta=0.58$.
 
 First three columns of table \ref{tb:importance} represent, respectively, indexes of nodes, vector $v$, and vector $\ell$ for the network
 in the figure \ref{fig:network}.
 Fourth and fifth columns represent, respectively, the values of $q_v^*$ and $q_{\ell}^*$ for a fixed budget $C=3.962$.
 These values highlight, also in this case an important difference between the two optimal solutions: notice that $q_v^*$ is greater
 than $1$ only for $6$ agents while $q_{\ell}^*$ is greater than 1 for all the agents.
 
 \begin{table}[!h]
  \begin{center}
  \caption{Optimal protections $q_v^*$ and $q_{\ell}^*$.}\label{tb:importance}
  \begin{tabular}{|l|l|l|l|l|} 
   \hline
   $i$ & $v$ & $\ell$ & $q_v^*$ & $q_{\ell}^*$ \\ 
   \hline\hline
    1  &  0.0575  &  0.4730  &  1.0000  &  1.0294\\
    2  &  0.1334  &  0.6965  &  1.6485  &  1.5159\\
    3  &  0.0901  &  0.4941  &  1.1137  &  1.0754\\
    4  &  0.0950  &  0.5428  &  1.1734  &  1.1813\\
    5  &  0.1124  &  0.5577  &  1.3891  &  1.2138\\
    6  &  0.0708  &  0.4632  &  1.0000  &  1.0081\\
    7  &  0.0792  &  0.5104  &  1.0000  &  1.1109\\
    8  &  0.0805  &  0.4775  &  1.0000  &  1.0393\\
    9  &  0.1066  &  0.5521  &  1.3178  &  1.2016\\
   10  &  0.1056  &  0.6749  &  1.3054  &  1.4689\\
   11  &  0.0688  &  0.5425  &  1.0000  &  1.1807\\
   \hline
  \end{tabular}
  \end{center}
 \end{table}

 Figure \ref{fig:q_functionC} shows entries of $q_v^*$ (in blue) and $q_{\ell}^*$ (in red) as a function of $C$ from $C=\sqrt{n}$ to $C=n$.
 There is a flat line until the high regime is not reached.

\begin{figure}[!ht]
 \begin{center}
 \includegraphics[width=8.4cm]{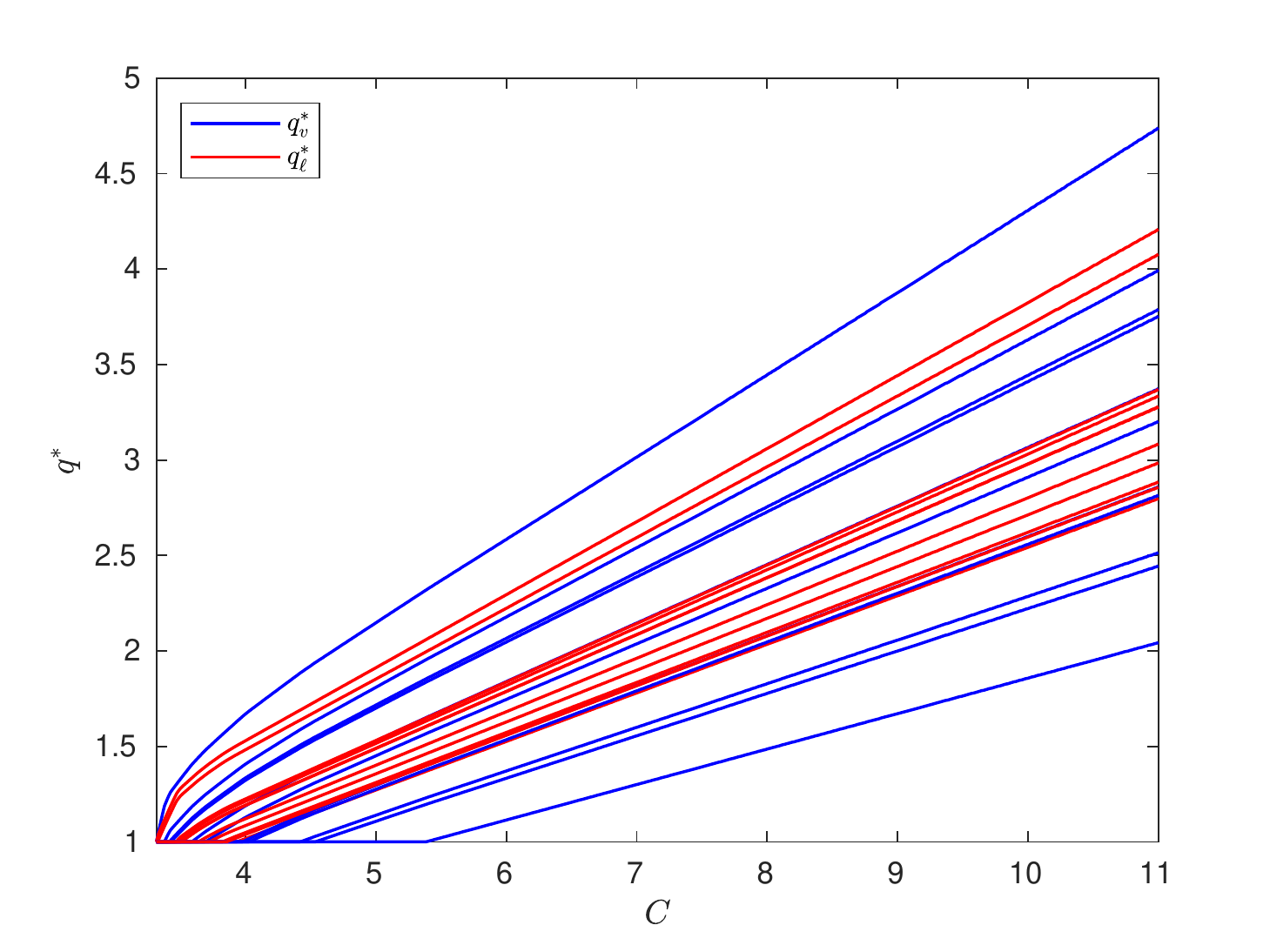}
 \caption{Entries of optimal protections $q_v^*$ and $q_{\ell}^*$ as a function of the budget $C$. Origin is located in $(\sqrt{n},1)$.}
 \label{fig:q_functionC}
 \end{center}
\end{figure}

\end{example}
\begin{example}
In this example we consider a network derived from a concrete dataset describing production interconnections among US firms published in 2002 by
the Bureau of Economic Analysis. 
This dataset has been used in (\cite{acemoglu2012network}) to which we refer for further description. The model consists
of an input-output matrix $A$ (Commodity-by-Commodity Direct Requirements Tables) that represents direct interactions between $n=417$ commodities.
The value $A_{ij}$ gives the input share of commodity $j$ used in the production of commodity $i$. $P$ is the normalized version of it that is used
in analyzing production networks with constant returns to scale technologies.

We have set $\beta=0.58$ taking this value from the aforementioned article.
 
For this example we have plot in figures \ref{fig:phiv_functionC} and \ref{fig:phil_functionC} the optimal value $\lambda(C)^*$ and we have compared
it against the protection extended to all nodes in a way proportional to their centrality. We indicate such value $\lambda(C)^{\rm diff}$.
Figures \ref{fig:phiv_ratio} and \ref{fig:phil_ratio} represent the ratio between the strategies.

\begin{figure}[!ht]
 \begin{center}
 \includegraphics[width=8.4cm]{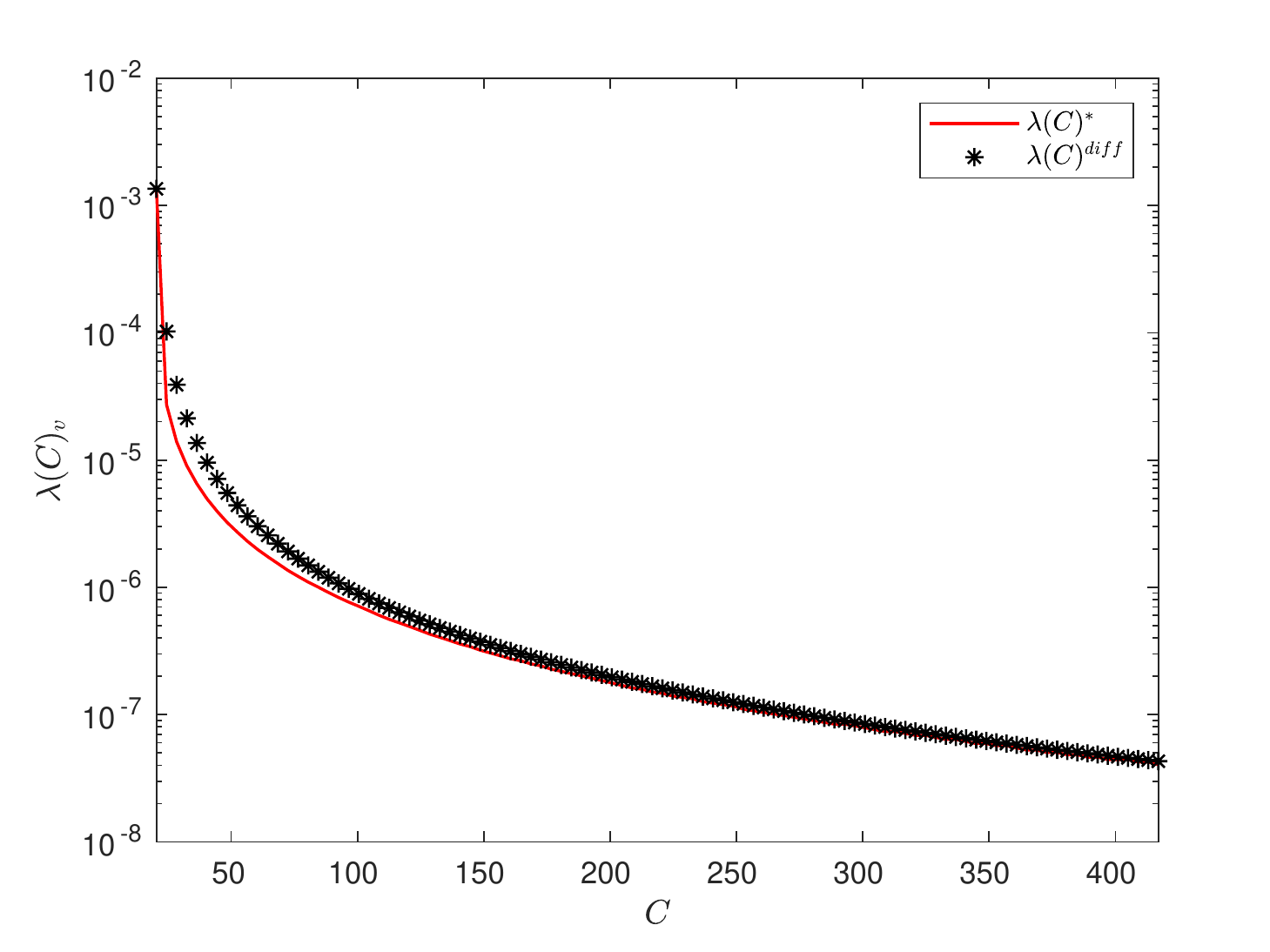}
 \caption{Optimal variance of the arithmetic mean $\lambda(C)_v^*$ compared with $\lambda(C)_v^{\rm diff}$.}
 \label{fig:phiv_functionC}
 \end{center}
\end{figure}

\begin{figure}[!ht]
 \begin{center}
 \includegraphics[width=8.4cm]{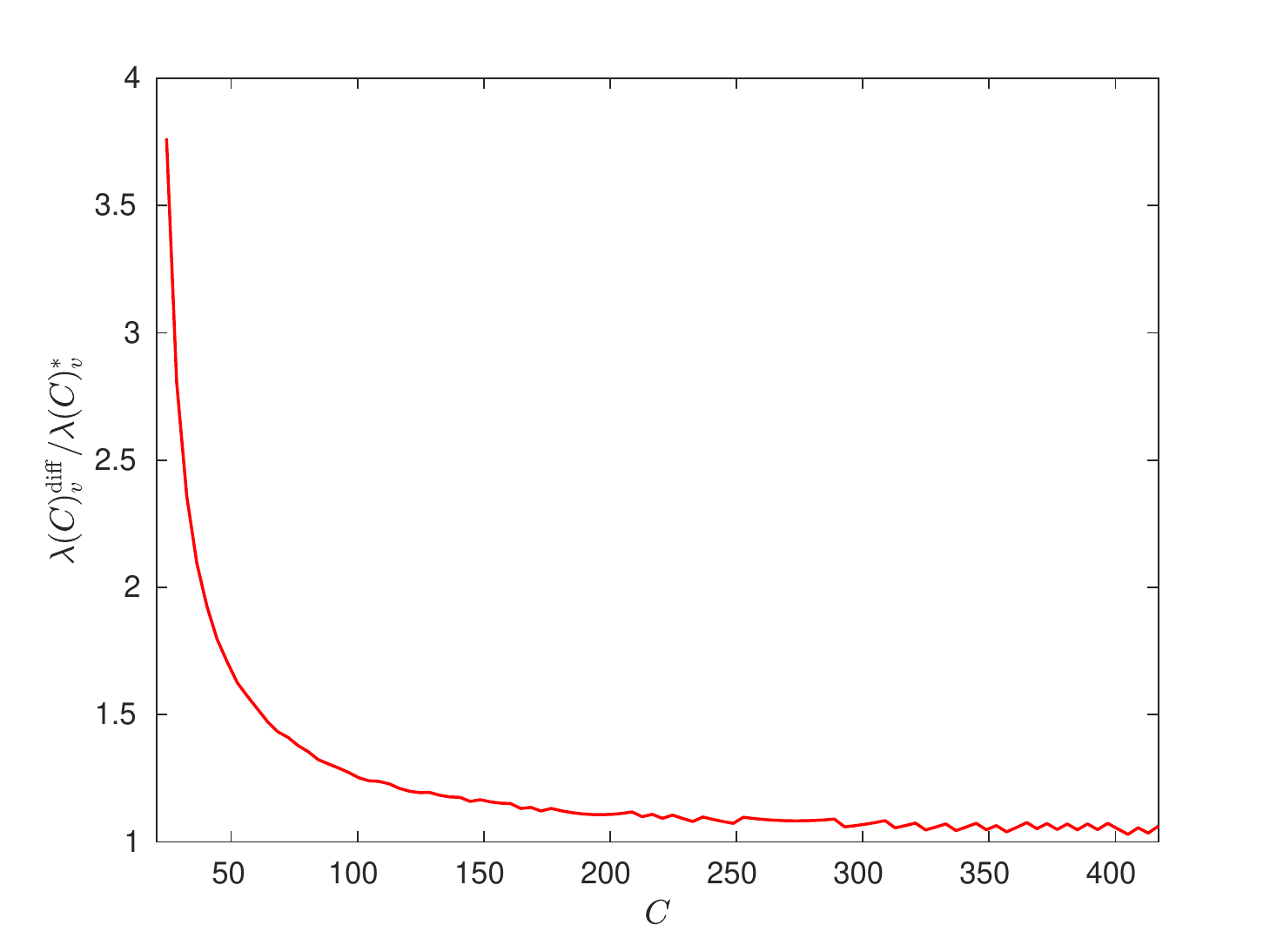}
 \caption{Ratio between $\lambda(C)_v^*$ and $\lambda(C)_v^{\rm diff}$ as function of $C$.}
 \label{fig:phiv_ratio}
 \end{center}
\end{figure}

\begin{figure}[!ht]
 \begin{center}
 \includegraphics[width=8.4cm]{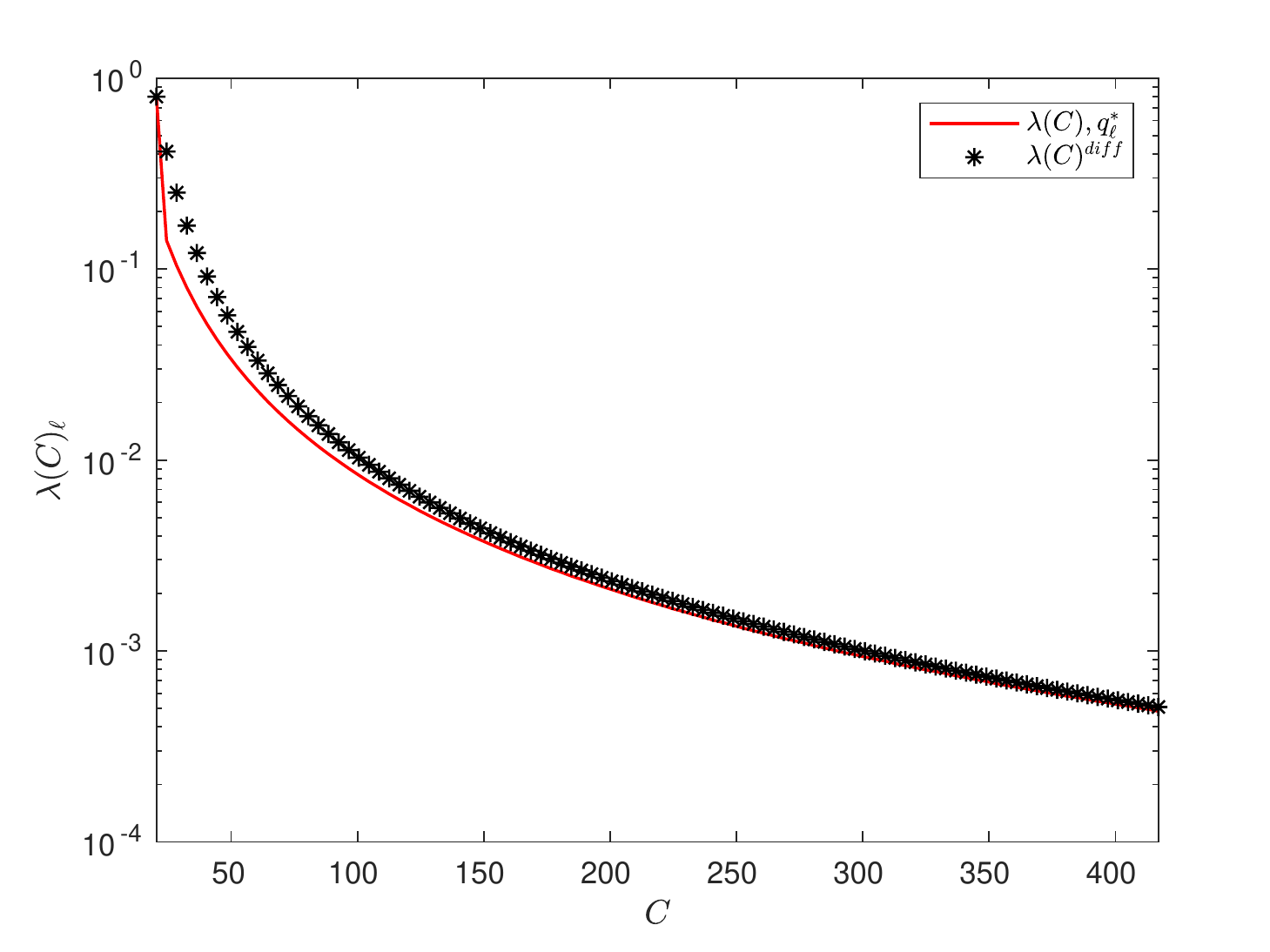}
 \caption{Optimal variance of $x$ $\lambda(C)_{\ell}^*$ compared with $\lambda(C)_{\ell}^{\rm diff}$.}
 \label{fig:phil_functionC}
 \end{center}
\end{figure}

\begin{figure}[!ht]
 \begin{center}
 \includegraphics[width=8.4cm]{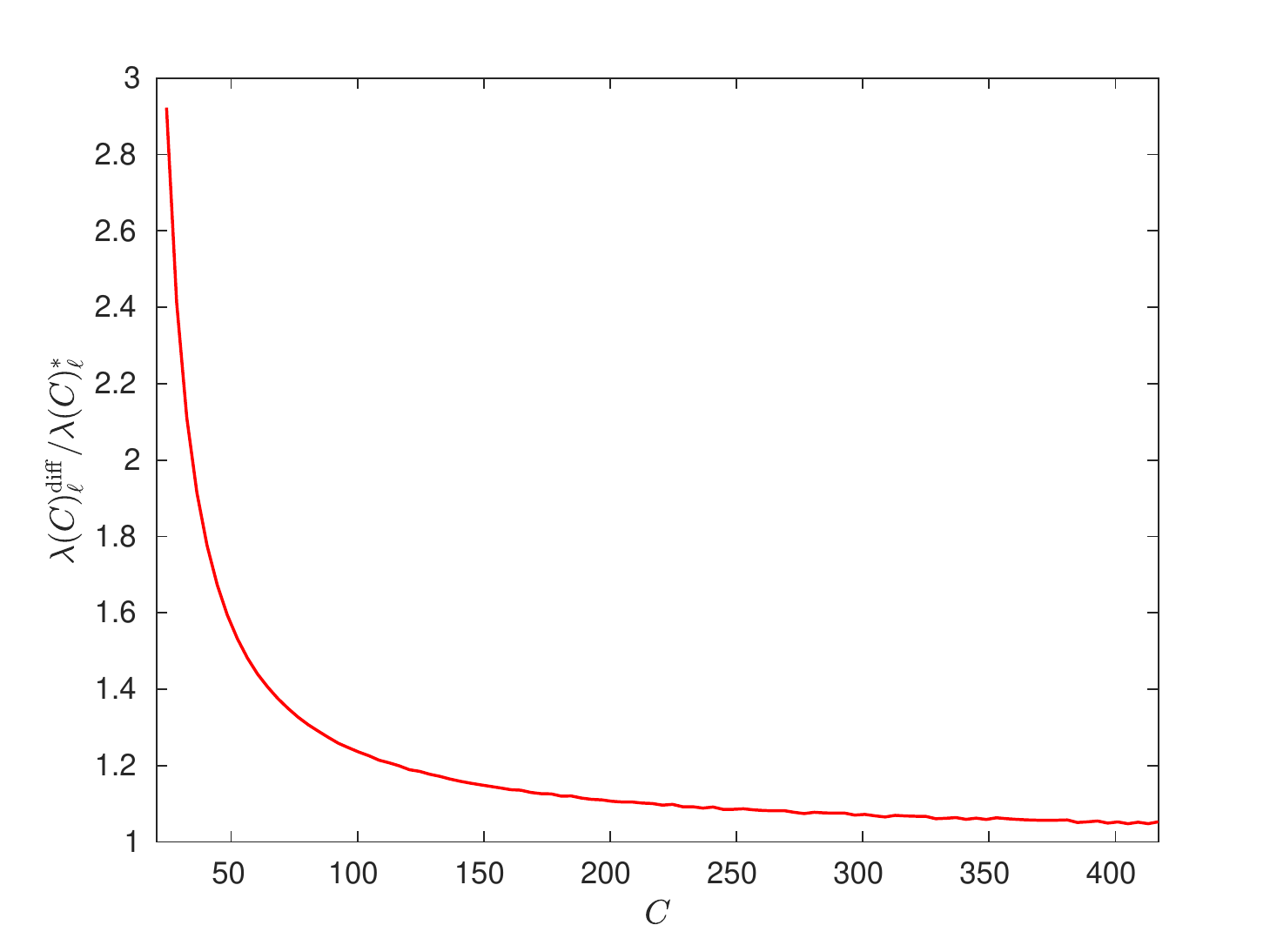}
 \caption{Ratio between $\lambda(C)_{\ell}^*$ and $\lambda(C)_{\ell}^{\rm diff}$ as function of $C$.}
 \label{fig:phil_ratio}
 \end{center}
\end{figure}

\end{example}

\section{Conclusion and future directions}\label{conclusion}

We have presented a new class of problems where shocks are complementary paired with protections and casted into
an adversarial min-max problem. We have explicitly solved the optimization problem, for the case of two possible aggregative performance functionals,
and studied the form of the optimal protection vector that exhibits a waterfiling shape that crucially depends on the aggregate performance chosen.

Future research steps will be in the direction of considering correlated shocks and more realistic constraints in the protection vector (e.g. heterogeneity of nodes).  

\bibliography{ifac_DamonteComoFagnani}

\end{document}